\documentstyle[amstex,12pt,bezier]{article}

\normalbaselines
\newenvironment{pf}%
    {\par\noindent{\it Proof\/: }\nopagebreak\normalsize}%
    {\hfill\linebreak[2]\hspace*{\fill}$\Box$\\[12pt]}

\newtheorem{thm}{Theorem}[section]
\newtheorem{prop}[thm]{Proposition}

\newtheorem{cor}[thm]{Corollary}
\newtheorem{lem}[thm]{Lemma}
\newtheorem{defn}[thm]{Definition}
\newtheorem{defns}[thm]{Definitions}
\newtheorem{nott}[thm]{Notation}
\newtheorem{rmk}[thm]{Remark}
\newtheorem{rmks}[thm]{Remarks}

\newcommand{\K}{{\cal K}}

\newcommand{\Z}{{\Bbb Z}}

\newcommand{\R}{{\Bbb R}}

\newcommand{\r}{{\cal R}}

\newcommand{\gtp}{{\frak p}}

\newcommand{\gtm}{{\frak m}}

\newcommand{\Sing}{{\rm Sing}}
\newcommand{\di}{{\rm dim}}

\newcommand{\spr}{{\rm Spec}_r}
\newcommand{\frz}{\partial _{\rm z}}
\newcommand{\adz}{{\rm Adh}_{\rm z}}

\newcommand{\ol}{\overline}

\newcommand{\notsubset}{\subset\kern-10pt/\;}
\newcommand{\notsupset}{\supset\kern-12pt/\;}
\newcommand{\om}{{\mit\Omega}}
\newcommand{\s}{\sigma}
\newcommand{\al}{\alpha}
\begin{document}
\title{AN ALGORITHMIC  CRITERION\\
  FOR BASICNESS IN DIMENSION 2}
\author{Francesca Acquistapace\thanks{Members of CNSAGA of CNR, partially
supported by MURST}\\
        Fabrizio Broglia$ ^*$ \\
       {\small Dipartimento di Matematica,}\\
       {\small Universit\`a di Pisa, I56100 Pisa, Italy}\\[12pt]
       M. Pilar V\'elez\thanks{Partially supported by DGICYT and a
Beca-Complutense}\\
       {\small Departamento de Geometr\'\i a y Topolog\'\i a, Facultad de
Matem\'aticas,}\\
       {\small Universidad Complutense de Madrid, 28040 Madrid, Spain}\\[12pt]}

\date{December,13,1993}
\maketitle
\vspace{12pt}

\centerline{\sl Mathematics Subject Classification: \rm 14P20}

\vspace{12pt}
\section*{Introduction.}

In this paper we give a constructive procedure to
check basicness of open (or closed) semialgebraic sets in a compact,
non singular, real algebraic surface $X$. It is rather clear that if a
semialgebraic set $S$ can be separated
from each connected component of $X\setminus(S\cup\frz S)$ (when $\frz S$
stands for the Zariski closure of $(\ol S\setminus{\rm Int}(S))\cap{\rm
Reg}(X)$), then $S$ is basic.
This leads to associate  to $S$ a finite family of sign distributions on
$X\setminus\frz S$;
we prove the equivalence between basicness and two properties of these
distributions, which can be tested by an algorithm described in
[ABF].

By this method we find for surfaces a general result of [AR2] about the
``ubiquity
of Lojasiewicz's example" of non basic semialgebraic sets
(2.9). There is a close relation between these two properties and the
behaviour of fans in the algebraic functions field of $X$ associated to a real
prime divisor (Lemmas 3.5 and 3.6). We use this fact to get an easy proof
(Theorem 3.8), for
a general surface $X$, of the well known 4-elements fan's criterion for
basicness (see [Br] and [AR1]). Furthermore, if the criterion fails,
using the description of fans in dimension 2 [Vz],
we find an algorithmic method to exhibit the failure. Finally, exploiting this
thecnics
of sign distribution we give one improvement of the 4-elements fan's criterion
in [Br]
to check if a semialgebraic set is principal.

In fact, one goal of this paper is to give purely geometric proofs, in the case
of surfaces, of the theory of fans currently used in Semialgebraic Geometry
([Br]).
 In particular, we only need the definitions of that theory. It is
at least remarkable that while the notion of fan is highly geometric in nature
all
known proofs of the main results are pure quadratic forms theory.

\section{Geometric review of basicness}

Let $X\subset \R^n$ be an algebraic surface. Denote by $\r(X)$  the ring of
regular
functions on $X$.
Let $S\subset X$ be a {\em semialgebraic set}, that is
$$S=\bigcup_{i=1}^p\{x\in X:f_{i1}(x)>0,\dots,f_{ir_i}(x)>0,g_i(x)=0\}$$
with $f_{i1},\dots,f_{ir_i},g\in \r(X)$, for $i=1,\dots,p$.\\
We will simply write: $S=\bigcup\{f_{i1}>0,\dots,f_{ir_i}>0,g_i=0\}$.

\begin{defn} A semialgebraic set $S$ is {\em basic open} (resp. {\em basic
closed})
if there exist $f_1,\dots,f_r\in \r(X)$ such that
$$S=\{x\in X:f_1(x)>0,\dots,f_r(x)>0\}$$
$$({\it resp.}\; S=\{x\in X:f_1(x)\geq 0,\dots,f_r(x)\geq 0\})$$
\end{defn}

\begin{defn} A semialgebraic set $S\subset X$ is {\em generically basic} if
there
exist a Zariski closed set $C\subset X$, with $\di(C)\leq 1$, such that
$S\setminus C$ is basic open.
\end{defn}

Denote by $\frz S$ the Zariski closure of the set $\partial(S)=(\ol
S\setminus{\rm Int}(S))\cap{\rm Reg}(X)$.
It is known that in dimension 2 {\em basic} and {\em generically basic} have
almost the same meaning (see [Br]). We give here a direct proof of this fact.

\begin{lem}
 Let $S$ be an open semialgebraic set in $X$, if $S$ is generically
basic then $S\cap \frz S$ is a finite  set.
\end{lem}

\begin{pf}
Suppose there exist $f_1,\dots,f_s\in \r(X)$ and an algebraic set $C\subset X$
such that $\di(C)\leq 1$ and
$$S\setminus C=\{f_1>0,\dots,f_s>0\}.$$
Suppose also that there is an irreducible component $H$ of $\frz S$ such that
$\di(H\cap S)=1$. Let $\gtp$ be the ideal of $H$ in $\r(X)$ and pick any
$x_0\in {\rm Reg}(H)$.
Then $\r(X)_{x_0}$ (the localization of $\r(X)$ at the maximal ideal $\gtm$ of
$x_0$, so $\gtp\subset\gtm$) is a factorial ring  and ${\rm
ht}(\gtp\r(X)_{x_0})=1$, hence
$\gtp\r(X)_{x_0}=h\r(X)_{x_0}$ for some $h\in \gtp$. Take
$g_1,\dots,g_r\in\r(X)$
such that $\gtp\r(X)=(g_1,\dots,g_r)$, then there exist $\lambda_i,s_i\in\r(X)$
with $s_i(x_0)\not= 0$ (in particular, $s_i\not\in\gtp$) such that
$s_ig_i=\lambda_ih$,
for $i=1,\dots,r$.\\
\indent Consider $U=X\setminus\{s_1\cdots s_r=0\}$ which is Zariski open in
$X$.
Then we have $\gtp\r(X)_x=(h)\r(X)_x$ for all $x\in U_1=U\cap{\rm Reg}(H)$.\\
\indent Take now $x_0\in U_1$, we have $f_j=\rho_jh^{\alpha_j}$ with
$\alpha_j\geq 0$
and $\rho_j\in \r(X)_{x_0}$ such that $h$ does not divide $\rho_j$ (in
particular, $\rho_j\not\in \gtp\r(X)_{x_0}$).
Then $\rho_j=p_j/q_j$, $j=1,\dots,s$, where $p_j,q_j\in\r(X)$ and
$p_j(x_0)\not=0,q_j(x_0)\not=0$.
Let $U_2$ be the Zariski open set $U_1\setminus\{q_1\cdots q_s=0\}\subset H$,
then for all $x\in U_2$,
$q_jf_j=p_jh^{\alpha_j}$, $p_j,q_j$
do not change sign in a neighbourhood of $x$, for $j=1,\dots,s$ and $h$
changes sign in a neighbourhood of $x$, since locally $h$ is a parameter of $H$
in $x$. Hence for any $j=1,...,s$, $f_j$ changes sign in a neighbourhood of $x$
if $\alpha_j$ is odd and does not change sign if $\alpha_j$ is even.\\
\indent Using the fact that $\di(S\cap H)=1$, $\di(C)\leq 1$ and
$s_i,p_j,q_j\not\in\gtp$, there exist a
Zariski dense open set $\om$ in $U_2$ such that $f_j$ does not change
sign through $\om$, for all $j=1,\dots,s$, then $\alpha_j$ is even for all
$j=1,\dots,s$.
But  also there is a Zariki dense open set $\om'$ in $U_2$
such that $\om'\subset \overline S\setminus S$, then there is $l\in
\{1,\dots,s\}$
such that $f_l$ changes sign through $\om'$, and $\alpha_l$ is odd, which is
impossible.
\end{pf}

\begin{prop}
Let $S$ be a semialgebraic set in $X$.\\
\indent (1) Let $S$ be open. Then $S$ is generically basic if and only if there
exist
$p_1,\dots,p_l\in \frz S$ such that $S\setminus\{p_1,\dots,p_l\}$ is basic
open.\\
\indent (2) $S$ is basic open if and only if $S$ is generically
basic and $S\cap\frz S=\emptyset$.\\
\indent (3) Let $S$ be closed. Then $S$ is generically basic if and only if $S$
is basic closed.
\end{prop}

\begin{pf}  First we prove {\it (1)}. The ``if part'' is trivial.
Suppose $S$ to be open and generically basic, that is there are
$f_1,...,f_r\in\r(X)$
and an algebraic set $C\subset X$ with $\di(C)\leq 1$ such that
$$S\setminus C=\{f_1>0,\dots,f_r>0\}.$$
\noindent We suppose first that $C$ is a curve and
we can also suppose  that $\{f_1\cdots f_r=0\}\subset C$.\\
\indent $C\cap S$ is an open semialgebraic set in $C$, then using [Rz, 2.2],
there exist $g_1\in \r(C)$ such that
$$C\cap S=\{x\in C:g_1(x)>0\}$$
$$\ol{C\cap S}=\{x\in C:g_1(x)\geq 0\}$$
\noindent choose $g_1$ to be the restriction of a regular function $g\in
\r(X)$.\\
\indent For each $i=1,\dots,r$ consider the open sets
$B_i=\{x\in X:f_i(x)<0\}$
\noindent and the closed sets
$T_i=(\ol S\cap\{g\leq 0\})\cup(\ol{B_i}\cap\{g\geq 0\}).$\\
\indent Applying [BCR, 7.7.10] to $T_i$, $f_i$ and $g$ for $i=1,...,r$,
we can find $p_i,q_i\in \r(X)$, with $p_i>0$, $q_i\geq 0$, such that\\
\indent {\em (i)} $F_i=p_if_i+q_ig$ has the same signs as $f_i$ on $T_i$;\\
\indent {\em (ii)} The zero set $Z(q_i)$ of $q_i$ verifies,
$Z(q_i)=\adz(Z(f_i)\cap T_i)$.\\
We remark the following:\\
\indent {\em a)} $F_i(S\setminus C)>0$ for $i=1,\dots,r$, since $F_i$ has the
same
signs as $f_i$ on $T_i\cap S$ and outside is the sum of a strictly positive
function and a nonnegative one.\\
\indent {\em b)} $F_i(B_i)<0$ for $i=1,\dots,r$, by the same reasons.\\
\indent {\em c)} $Z(q_i)\subset \frz S$. In fact, denote
$$Z_1^i=Z(f_i)\cap \ol S\cap\{g\leq 0\}$$
$$Z_2^i=Z(f_i)\cap \ol{B_i}\cap\{g\geq 0\}$$
\noindent then we have
$Z(q_i)=\adz(Z_1^i)\cup\adz(Z_2^i)\subset Z(f_i)\subset C.$
 Since $g$ is positive on $C\cap S$ and $Z_1^i\subset C\cap \{g\leq 0\}$,
we have $Z_1^i\cap S=\emptyset$, hence $Z_1^i\subset \partial(S)$. Indeed,
since
$B_i\cap S=\emptyset$ and $S$ is open we have
$\ol{B_i}\cap S=\emptyset$. Moreover $C\cap \ol{B_i}\subset\{g\leq 0\}$,
hence $Z_2^i\subset\{g=0\}\cap C\subset \ol S$, then
$Z_2^i\subset\partial(S)$.\\
\indent From these remarks, denoting $Z=\bigcup_{i=1}^rZ(q_i)$, we have
$$S\setminus Z=\{F_1>0,\dots,F_r>0\}.$$
\noindent In fact, if $x\in S\setminus C$ then $F_i(x)>0$ for $i=1,...,r$,
if $x\in (C\cap S)\setminus Z$ then $f_i(x)\geq 0$, $q_i(x)\geq 0$
and $g(x)>0$, hence $F_i(x)>0$ for $i=1,...,r$ and $x\in S\setminus Z$.
Otherwise, suppose $x\not\in S\setminus Z$, then $x\in (X\setminus S)\cup(S\cap
Z)$:
if $x\in X\setminus S$ there is $l\in \{1,\dots,r\}$ such that $f_l(x)\leq 0$,
and we can have $x\not\in C$ or $x\in C$,
in the first case $f_i(x)\not= 0$ for all $i$, so $x\in B_l$ and $F_l(x)<0$, in
the second case $g(x)\leq 0$, $q_l(x)\geq 0$, so
$F_l(x)\leq 0$; if $x\in S\cap Z$ there is $l\in \{1,\dots,r\}$ such that
$q_l(x)=0$, then
$f_l(x)=0$ and $F_l(x)=0$. In any way, there is $l$ such that $F_l(X)\leq 0$ if
$x\not\in S\setminus Z$.\\
\indent By 1.3 and  remark {\em c)} above we have
that there exist $p_1,\dots,p_l\in\frz S$
such that $\bigcup_{i=1}^rZ(q_i)\cap S=\{p_1,\dots,p_l\}$, hence
$$S\setminus\{p_1,\dots,p_l\}=\{F_1>0,\dots,F_r>0\}.$$
\indent If $C=\{a_1,\dots,a_m\}$ is a finite set we have to check that we can
throw out from $C$ all the $a_i$ which do not lie in $\frz S$. This can be done
as before by taking the function 1 at the place of $g$ and putting
$T_i=\ol{B_i}$.\\
\indent From {\it (1)} we have immediately {\it (2)}, because if $S$ is
generically
basic and $\frz S\cap S=\emptyset$, $S$ is basic open, since following the
proof above we have $Z(q_i)\cap S=\emptyset$ for $i=1,...,r$. On the countrary,
if $S$ is basic open then it is generically basic and $\frz S\cap S=\emptyset$
(because $\frz S\subset\{f_1\cdots f_r=0\}$ if $S=\{f_1>0,\dots,f_r>0\}$).\\
\indent Finally we prove {\it (3)}. The ``if part" is trivial. Then  suppose
$S$
to be closed and generically basic, i. e.
$$S\setminus C=\{f_1>0,\dots,f_r>0\}$$
\noindent with $f_1,\dots,f_r\in \r(X)$ and $\di(C)<2$. We can suppose
$\{f_1\cdots f_r=0\}\subset C$.\\
\indent $C\cap S$ is a closed semialgebraic in $C$, by [Rz, 2.2], there is
$g_1\in \r(C)$ such that
$$C\cap S=\{x\in C:g_1(x)\geq 0\}.$$
\noindent Take $g\in \r(X)$ as above, $f\in \r(X)$ a positive equation of $C$
and $T=S\cap \{g\leq 0\}$. Applying again [BCR, 7.7.10] to $T$, $f$ and $g$,
we find $p,q\in\r(X)$ with $p>0$, $q\geq 0$ such that\\
\indent {\em (i)} $h=pf+qg$ has the same sign as $f$ on $T$;\\
\indent {\em (ii)} $Z(q)=\adz(Z(f)\cap T)$.\\
Notice that $h(S)\geq 0$ and $Z(f)\cap T=C\cap S\cap\{g\leq
0\}\subset\{g=0\}\cap C$,
because $C\cap S\subset\{g\geq 0\}$; but $\{g=0\}\cap C$ is a finite set
contained in $S$, then $Z(q)$ is a finite set contained in $S$.\\
\indent We will prove that
$$S=\{f_1\geq 0,\dots,f_r\geq 0,h\geq 0\}.$$
\noindent In fact, if $x\in S\setminus C$ then $f_i(x)>0$ for all $i$ and
$h(x)>0$, if $x\in S\cap C$ then $f_i(x)\geq 0$ for all $i$ (because
$\di(C)<2$) and $h(x)\geq 0$,
hence $S\subset\{f_1\geq 0,\dots,f_r\geq 0,h\geq 0\}$. Otherwise suppose
$x\not\in S\supset S\setminus C$,
then there is $l\in \{1,...,r\}$ such that $f_l(x)\leq 0$,
if $x\not\in C$, $f_i(x)\not= 0$ for all $i=1,...,r$, then $f_l(x)<0$;
if $x\in C\setminus (C\cap S)$, $f_l(x)\leq 0$, $q(x)\not= 0$ and $g(x)<0$,
then $h(x)<0$.
\end{pf}

\begin{rmks}
{\em (1)} {\rm Let $S$ be a closed semialgebraic set, then $S$ is basic closed
if and only if $S\setminus\frz S$ is basic open.\\
\indent In fact, if $S$ is basic closed then it is generically basic, hence
$S\setminus\frz S$
is generically basic and as $(S\setminus\frz S)\cap\frz(S\setminus\frz
S)=\emptyset$,
$S\setminus\frz S$ is basic; on the contrary, if $S\setminus\frz S$ is basic
open then
$S$ is basic closed, since $S$ is closed.}

\indent {\em (2)} {\rm Let $S$ be a semialgebraic set, and let $S^\ast$ denote
the set
${\rm Int}(\ol S)$. Then, $S$ is basic open if and only if $S^\ast$ is
generically
basic and $S\cap\frz S=\emptyset$.\\
\indent In fact, $S$ and $S^\ast$ are generically equal.}
\end{rmks}

\section{Basicness and sign distributions}

We recall some definitions and results from [ABF]. Let
$X$ be a compact, non singular, real algebraic surface and $Y\subset X$ an
algebraic curve.\\
\indent Consider a {\em partial sign distribution} $\s$ on $X\setminus Y$,
which gives the sign $1$ to some connected components of $X\setminus Y$ (whose
union is denoted by $\s^{-1}(1)$) and the sign $-1$ to some others (whose union
is denoted
by $\s^{-1}(-1)$). So $\s^{-1}(1)$ and $\s^{-1}(-1)$ are disjoint open
semialgebraic
sets in $X$.

\begin{defns}
{\em (1)} A sign distribution $\s$ is {\em completable} if $\s^{-1}(1)$
and $\s^{-1}(-1)$ can be separated by a regular function, i.e. there is
$f\in\r(X)$
such that $f(\s^{-1}(1))>0$, $f(\s^{-1}(-1))<0$ and $f^{-1}(0)\supset Y$.
Briefly,
we say that {\em $f$ induces $\s$}.

{\em (2)} A sign distribution $\s$ is {\em locally completable} at a
point $p\in Y$ if there is $f\in\r(X)$ such that $f$ induces $\s$ on a
neighbourhood of $p$.

{\em (3)} An irreducible component $Z$ of $Y$ is a {\em type changing
component}
with respect to $\s$ if there exist two nonempty open sets $\om_1,\om_2\subset
Z\cap{\rm Reg}(Y)$ such that\\
\hspace{.2in} (a) $\om_1\subset\ol{\s^{-1}(1)}\cap\ol{\s^{-1}(-1)}$,\\
\hspace{.2in} ($b_+$) $\om_2\subset {\rm Int}(\ol{\s^{-1}(1)})$ \hspace{.2in}
or \hspace{.2in}
($b_-$) $\om_2\subset {\rm Int}(\ol{\s^{-1}(-1)})$.\\
If ($b_+$) (resp. ($b_-$)) holds we say that $Z$ is {\em positive type
changing}
(resp. {\em negative type changing}) with respect to $\s$.

{\em (4)} An irreducible component $Z$ of $Y$ is a {\em change component} if
there exist a nonempty open set $\om\subset Z\cap{\rm Reg}(Y)$ verifying (a).
\end{defns}

Completable sign distributions are characterized by the following theorem.

\begin{thm}
{\rm (See [ABF, 1.4 and 1.7])} Denote by $Y^c$ the union of the change
components of $Y$ with respect to $\s$. Then $\s$ is completable if and
only if\\
\indent {\em (1)} $Y$ does not have type changing component with respect to
$\s$;\\
\indent {\em (2)} $\s$ is locally completable at any point $p\in\Sing(Y)$;\\
\indent {\em (3)} There exist an algebraic curve $Z\subset X$ such that
$Z\cap(\s^{-1}(1)\cup\s^{-1}(-1))=\emptyset$ and $[Z]=[Y^c]$ in ${\rm
H}_1(X,\Z_2)$.\\
\indent Moreover condition {\em (2)} becomes condition {\em (1)} after the
blowings-up
of the canonical desingularization of $Y$, namely each point where $\s$ is not
locally completable corresponds to at least one type changing component of the
exceptional divisor with respect to the lifted sign distribution.
\end{thm}

\begin{prop}
{(\rm See [ABF, ``the procedure"])} Condition {\em (2)} of theorem 2.2 can be
tested,
without performing the blowings-up, by an algorithm that only  uses the Puiseux
expansions of the branches of $Y$ at its singular points.
\end{prop}

In fact, there are two algorithms:

\noindent \fbox{A1} (see [ABF, 2.4.19]) Given a branch $C$ of an algebraic
curve through a point
$p_0$ and an integer $\rho>0$, it is possible to find explicitely, in terms
of the Puiseux expansion of $C$, the irreducible Puiseux parametrizations of
all analytic arcs $\gamma$ through $p_0$ with the following properties:\\
\indent a) Denoting respectively by $\gamma_{\rho}$ and by $C_{\rho}$ the
strict transform
of $\gamma$ and $C$ after $\rho$ blowings-up in the standard resolution
of $C$ and by $D_{\rho}$ the exceptional divisor arising at the last
blowing-up,
 then $\gamma_{\rho-1}$ is parametrized by
\[ \left\{ \begin{array}{l}
x=t\\
y=at+\cdots\end{array}\right. \]
\noindent and $C_{\rho-1}\cap\gamma_{\rho-1}=0\in D_{\rho-1}$.\\
\indent b) $\gamma_{\rho-1}$ and $C_{\rho-1}$ have distinct tangents at 0.

\noindent \fbox{A2} (see [ABF, solution to problem 2]) Given an analytic arc
$\gamma$ and a region of $\R^2$ bounded by two analytic arcs
$\gamma_1,\gamma_2$,
with $\gamma,\gamma_1,\gamma_2$ through $0\in \R^2$, it is possible to decide,
looking at the Puiseux parametrizations, whether $\gamma$ crosses the region
or not.

Using \fbox{A1} and \fbox{A2} one can decide whether $D_{\rho}$ is positive
or negative type changing with respect to the lifted sign distribution
$\s_{\rho}$
without performing the blowings up, because for the family of arcs given by
\fbox{A1}, whose strict transform are transversal to $D_{\rho}$, we can
decide by \fbox{A2} whether or not $\s$ changes sign along some elements of the
family
and whether or not $\s$ has constant positive or constant negative sign along
some other ones.\\

Now let $S$ be an open semialgebraic set.

\begin{nott}
{\rm Let $A_1,\dots,A_t$ be the connected component of $X\setminus(S\cup\frz
S)$.
For each $i=1,\dots,t$ we denote by $\s_i^S$ (or simply $\s_i$ when there
is not risk of confusion) the following sign distribution on $X\setminus\frz
S$:
\begin{eqnarray*}
(\s_i^S)^{-1}(1)&=&S\setminus\frz S\\
(\s_i^S)^{-1}(-1)&=&A_i
\end{eqnarray*}    }
\end{nott}

\begin{lem}
Let $S$ be a semialgebraic set and $S^\ast={\rm Int}(\ol S)$. Then,
$\di(S^\ast\cap\frz S^\ast)=1$ if and only if there exists $i\in\{1,\dots,t\}$
such that $\frz S$ has a positive type changing component with respect to
$\s_i^S$.
\end{lem}

\begin{pf}
Suppose $\di(S^\ast\cap\frz S^\ast)=1$, then there is an irreducible component
$H$ of $\frz S^\ast$ such that $\di(H\cap S^\ast)=1$. So we can find an
open 1-dimensional set $\om\subset H\cap S^\ast$, hence
$$\om\subset{\rm Int}(\ol S)={\rm Int}(\ol{S\setminus\frz S})={\rm
Int}(\ol{\s_i^{-1}(1)})$$
\noindent for each $i=1,\dots,t$. And we can find another open 1-dimensional
set $\om'\subset\partial(S^\ast)\cap H$ such that
$$\om'\subset\ol{S^\ast}=\ol S=\ol{S\setminus\frz S}=\ol{\s_i^{-1}(1)}$$
\noindent for each $i=1,\dots,t$ and
$$\om'\subset X\setminus S^\ast=\bigcup_{i=1}^t\ol{A_i}\; .$$
\noindent Then $\om'\subset\bigcup_{i=1}^t(\ol{A_i}\setminus A_i)$, since
$A_i\cap S=\emptyset$
and $S$ is open. Hence, there exist a 1-dimensional open set
$\om''\subset\om'$ and $i_0\in\{1,...,t\}$ such that
$\om''\subset\ol{A_{i_0}}\setminus A_{i_0}$,
because $\di(\om')=1$ and we have a finite number of $A_i$.
Hence $\om''\subset\ol{\s_{i_0}^{-1}(1)}\cap\ol{\s_{i_0}^{-1}(-1)}$.
But since $H$ is a 1-dimensional component of $\frz S^\ast$ and $\frz S^\ast
\subset\frz S$, $H$ is an irreducible component of $\frz S$ of dimension 1;
then if we take $\om_1=\om''\cap{\rm Reg}(\frz S)$ and
$\om_2=\om\cap{\rm Reg}(\frz S)$, $H$ is a positive type changing
component with respect to $\s_{i_o}$.\\
\indent On the contrary suppose that $H$ is an irreducible component of
$\frz S$ which is a positive type changing component with respect to some
$\s_l$ ($l=1,...,t$). So there exist open sets $\om_1,\om_2\subset H$ of
${\rm Reg}(\frz S)$ such that
$\om_1\subset\ol{\s_l^{-1}(1)}\cap\ol{\s_l^{-1}(-1)}$
and $\om_2\subset{\rm Int}(\ol{\s_l^{-1}(1)})$. Then $\om_2\subset S^\ast$
and $\di(S^\ast\cap H)=1$. But $\om_1\subset\ol{S\setminus\frz S}=\ol {S^\ast}$
and $\om_1\subset\ol{A_i}$, then $\om_1\subset\ol{S^\ast}\setminus S^\ast$
because $X\setminus S^\ast=\bigcup\ol{A_i}$. So $H$ is an irreducible
component of $\frz S^\ast$, hence $\di(S^\ast\cap\frz S^\ast)=1$.
\end{pf}

\begin{prop}
Let $S$ be a semialgebraic set in the sphere ${\Bbb S}^2$ such that
$\frz S\cap S=\emptyset$ (resp. a closed semialgebraic set in ${\Bbb S}^2$).
Then $S$ is basic open (resp. closed) if and only if for each $i=1,\dots,t$,
$\s_i^S$ is
completable.
\end{prop}

\begin{pf}
It suffices to prove the result for a semialgebraic $S$ such that $\frz S\cap
S=\emptyset$,
because if $S$ is closed, applying 1.5 we have done.\\
\indent Suppose $S$ to be basic open, then $S^\ast$ is generically basic and
$\di(S^\ast
\cap\frz S^\ast)<1$. By lemma 2.5 no irreducible component of $\frz S$ is
positive type changing with respect to $\s_i^S$ for each $i=1,...,t$. But also
they cannot be negative type changing, because any curve in ${\Bbb S}^2$
divides
${\Bbb S}^2$ into connected components in such a way that none of them lies on
both sides of a branch of the curve (because the curves in ${\Bbb S}^2$ have
orientable neighbourdhoods).\\
\indent If for some $i\in\{1,...,t\}$, $\s_i^S$ were not locally completable
at some point $p\in\frz S$, we could find (Theorem 2.2) a non singular
surface $V$ together with a contraction $\pi:V\to{\Bbb S}^2$ of an algebraic
curve $E\subset V$ to the point $p$, with $\pi^{-1}(\frz S)$ normal crossing
in $V$, such that an irreducible component $D$ of $E$ would be type changing
with respect to $\s_i'=\s_i^S\cdot \pi$. But $p\not\in S\cup A_i$, so $\s_i'$
is defined on $V\setminus\pi^{-1}(\frz S)$ by
\begin{eqnarray*}
(\s_i')^{-1}(1)=&\pi^{-1}(S)=&T\\
(\s_i')^{-1}(-1)=&\pi^{-1}(A_i)&;
\end{eqnarray*}
\noindent and as $\pi:V\setminus(\pi^{-1}(\frz S))\to {\Bbb S}^2\setminus\frz
S$
is a biregular isomorphism, $\s_i'=\s_i^T$. So $T$ and $T^\ast$ are generically
basic, being biregularly
isomorphic to $S$, hence  $\di(\frz T^\ast\cap
T^\ast)<1$; so $D$ cannot be positive type changing by lemma 2.5. We have to
exclude that $D$ is negative type changing. Suppose it is so; then
we can find two open sets $\om,\om'\subset D$ such that $\pi^{-1}(A_i)$ lies on
both sides of $\om$ and
$\om'$ divides $T$ from $\pi^{-1}(A_i)$. Then there exists an irreducible
component $Z$ of $\frz S$ such that its strict transform $Z'$ crosses $D$
between $\om$ and $\om'$.
Then $Z'$ must cross $\pi^{-1}(A_i)$, because this open set is connected. This
is impossible since $A_i$
does not lie on both sides of $Z$. Then no irreducible component of $\frz S$ is
type changing with
respect to $\s_i^S$ and $\s_i^S$ is locally completable at any $p\in\frz S$
for all $i$; so $\s_i^S$ is completable, since ${\rm H}_1({\Bbb S}^2,\Z_2)=0$
(2.2).\\
\indent Suppose now that $\s_i^S$ is completable for each $i=1,...,t$ and
let $f_i\in\r({\Bbb S}^2)$ be a regular function inducing $\s_i^S$. Clearly
$S\subset\{f_1>0,\dots,f_t>0\}.$
But if $x\not\in S$ then $x\in A_l$ for some $l=1,...,t$ or $x\in\frz S$.
If $x\in A_l$, then $f_l<0$; if $x\in\frz S$, then $f_i(x)=0$ for all
$i=1,...,t$,
since $f_i$ vanishes on $\frz S$ by the very definition of completability.
So $S=\{f_1>0,\dots,f_t>0\}$ and it is basic.
\end{pf}

Before proving an analogous result for a general surface we need a lemma.

\begin{lem}
Let $S\subset X$ be an open semialgebraic set such that $S=S^\ast$,
$\frz S\cap S=\emptyset$ and $\frz S$ is normal crossing. Let $H$ be an
irreducible component of $\frz S$. Then there exist a non singular algebraic
set $H'\subset X$ such that\\
\indent (a) $[H]=[H']$ in ${\rm H}_1(X,\Z_2)$,\\
\indent (b) $H$ and $H'$ are transversals,\\
\indent (c) $H'\cap S=\emptyset$.
\end{lem}

\begin{pf}
$H$ is a smooth compact algebraic curve, composed by several ovals and locally
disconnects its neighbourhood; moreover $S$ is locally only on one side of $H$,
because $S=S^\ast$ and $H\cap S^\ast=\emptyset$, at least outside a
neighbourhood
of the singular points lying in $H$, namely the points where $H$ crosses an
other component of $\frz S$. At each of these points we have only two possible
situations (see figures 2.7).
 From this it is clear how to construct a smooth differentiable curve $C_H$
with $[C_H]=[H]$ in ${\rm H}_1(X,\Z_2)$ such that $C_H$ is transversal to each
irreducible component of
$\frz S$ and $C_H\cap S=\emptyset$, using a suitable tubular neighbourhood of
$H$.

\begin{center} \setlength{\unitlength}{1mm} \begin{picture}(100,40)
\put(0,20){\line(1,0){40}}
\put(60,20){\line(1,0){40}}
\put(20,0){\line(0,1){40}}
\put(80,0){\line(0,1){40}}
\put(78,4){\makebox(0,0){$ _H$}}
\put(18,4){\makebox(0,0){$ _H$}}
\put(5,35){\makebox(0,0){$S$}}
\put(65,35){\makebox(0,0){$S$}}
\multiput(0,22)(0,2){5}{\line(1,0){19.5}}
\multiput(10,32)(0,2){4}{\line(1,0){9.5}}
\multiput(60,22)(0,2){5}{\line(1,0){19.5}}
\multiput(70,32)(0,2){4}{\line(1,0){9.5}}
\multiput(80.5,2)(0,2){9}{\line(1,0){19.5}}
\put(23,2){\line(0,1){36}}
\put(26,35){\makebox(0,0){$ _{C_H}$}}
\put(89,35){\makebox(0,0){$ _{C_H}$}}
\bezier{100}(85,39)(85,25)(80,20)
\bezier{100}(80,20)(75,15)(75,1)
\put(20,-4){\makebox(0,0){{\small Figure 2.7.a}}}
\put(80,-4){\makebox(0,0){{\small Figure 2.7.b}}}
\end{picture}
\end{center}
\vskip .2cm

 \indent Now we want to approximate $C_H$ by a nonsingular algebraic curve $H'$
with the same properties.
 To do this we use the fact that $[H]$ gives a strongly algebraic line
bundle $\pi:E\rightarrow X$ on $X$ (see [BCR, 12.2.5]). For this line bundle,
$H$ is the zero set of an algebraic section $h$ and $C_H$ is the zero set
of a ${\cal C}^\infty$ section $c$. Let $Q_1,\dots,Q_k$ be the points in the
set
$H\cap C_H$. We can take a finite open covering $V_1,\dots,V_l$ of $X$ with
the following properties:\\
\indent 1) $Q_i\in V_i$ for $i=1,...,k$ and $V_1,\dots,V_k$ are pairwise
disjoint.\\
\indent 2) For each $j=1,...,l$ there exist an algebraic section $s_j$ of $E$
such that $s_j(x)\not=0$ for each $x\in V_j$ ($s_j$ generates $E_x$ for each
$x\in V_j$).\\
\indent Take a ${\cal C}^\infty$ partition of the unity
$\{\varphi_1,\dots,\varphi_l\}$
associated to the covering, with the property that for $i=1,...,k$,
$\varphi_i^{-1}(1)$
is a closed neighbourhood of $Q_i$ in $V_i$, $\varphi_i(Q_j)=0$ for $j\not= i$.
For $j=1,...,l$ we can write
$c\,\vline_{V_j}=\alpha_j s_j,$
\noindent with $\alpha_j\in {\cal C}^\infty(V_j)$ and $\alpha_j(Q_j)=0$ if
$j=1,...,k$, because $s_j$ generates the fiber. Then we have
$$c=\sum_{j=1}^l\varphi_jc=\sum_{j=1}^l(\varphi_j\alpha_j)s_j$$
\noindent and the smooth functions $\beta_j=\varphi_j\alpha_j\in{\cal
C}^\infty(X)$
vanishes at $Q_i$ for $i=1,...,k$.\\
\indent By a classical relative approximation theorem (see [BCR, 12.5.5]) we
can approximate $\beta_j$ on $X$ by a regular function $f_j$ on $X$ vanishing
at $\{Q_1,\dots,Q_k\}$. Then the algebraic section
$$s=\sum_{j=1}^lf_js_j$$
\noindent has an algebraic zero set $H'$ passing through $Q_1,\dots,Q_k$.
Moreover $H'\cap S=\emptyset$, because $H'$ is very close to $C_H$, and
$[H']=[H]$ in ${\rm H}_1(X,\Z_2)$.
\end{pf}

Finally we have.

\begin{thm}
Let $X$ be a compact, non-singular, real algebraic surface and $S\subset X$
be a semialgebraic set with $\frz S\cap S=\emptyset$ (resp. $S\subset X$
be a closed semialgebraic set).\\
\indent Then $S$ is basic if and only if for each $i=1,\dots,t$ the sign
distribution $\s_i^S$ verifies the following two properties:\\
\indent (a) No irreducible component of $\frz S$ is positive type changing
with respect to $\s_i^S$.\\
\indent (b) No irreducible component of the exceptional divisor of a
standard resolution \linebreak $\pi:V\to X$ of $\frz S$ is positive type
changing
with respect to $\s_i'=\s_i^S\cdot\pi$.
\end{thm}

\begin{pf}
It suffices to prove for $S$ with $S\cap\frz S=\emptyset$, because this implies
that $S$ is open and by remark 1.5.1 we have done for $S$ closed.\\
\indent The ``only if part" is the same as for ${\Bbb S}^2$, without the
argument
proving there are not negative type changing component. For the
 ``if part" we can reason as follows.\\
\indent Let $\pi:V\to X$ be the standard resolution of $\frz S$. Denote by $Y$
the curve $\pi^{-1}(\frz S)$ (see [EC] or [BK]), $Y$ is normal crossing in $V$.
Each irreducible component of $Y$ is a non-singular curve consisting possibly
of several ovals. Each of this ovals can have an orientable neighbourhood (in
which case it is homologically trivial) or a non orientable neighbourhood
isomorphic to the M\" obius band.\\
\indent Denote as it is usual $\s_i'$ the sign distribution $\s_i^S\cdot \pi$
and consider the sets $T=\pi^{-1}(S)$ and $T^\ast={\rm Int}(\ol T)$. Conditions
{\em (a)} and {\em (b)} for $\s_i^S$ imply that $Y$ has not positive type
changing components with respect to $\s_i'$. So by 2.5, $\frz T^\ast\cap
T^\ast$ is a finite
set, and as $\frz T^\ast\subset Y$ has not isolated points and $T^\ast$ is
open, we have $T^\ast\cap \frz T^\ast=\emptyset$.\\
\indent Consider now the sign distributions $\s_j^{T^\ast}$ on $V\setminus
\frz T^\ast$ defined as before for $j=1,\dots,l$, where $l$ is the number
of connected component of $V\setminus(T^\ast\cup\frz T^\ast)$; clearly they do
not have positive type changing components. Apply 2.7 to each irreducible
component $H$ of $\frz T^\ast$ being
a change component for some $\s_j^{T^\ast}$. Then, we find a non-singular
algebraic set $H'$ such that $H'\cap T^\ast=\emptyset$ and $[H]=[H']$ in
${\rm H}_1(V,\Z_2)$. The union of all $H'$ gives an algebraic set $Z\subset
V$.\\
\indent Remark that if $H$ is a negative type changing component with respect
some $\s_j^{T^\ast}$ then this phenomenon occurs along an oval of $H$ whose
neighbourhood is a M\" obius band, because in the other case the two sides of
the
oval whould be in different connected components of $V\setminus\frz T^\ast$.
Take the sign distributions $\tau_K$ on $V\setminus(\frz T^\ast\cup Z)$,
$k=1,\dots,m$,
defined by
\begin{eqnarray*}
(\tau_k)^{-1}(1)&=&T^\ast\\
(\tau_k)^{-1}(-1)&=&B_k,
\end{eqnarray*}
\noindent where $B_1,\dots,B_m$ are the connected components of $V\setminus
(T^\ast\cup\frz T^\ast\cup Z)$.\\
\indent We claim that $\tau_k$ is completable for each $k=1,...,m$. In fact,
we prove that conditions (1), (2) and (3) of 2.2 are verified.\\
\indent (1) No irreducible component of $\frz T^\ast\cup Z$ is neither
positive nor negative type changing with respect $\tau_k$: this is true because
 now the sign $-1$ can occur at most on one side of each irreducible component
$H$ of $\frz T^\ast\cup Z$ (a M\" obius band is divided by two transversal
generators of its not vanishing homological class into two connected
components), then
there are not negative type changing components with respect $\tau_k$.\\
\indent (2) $\tau_k$ is completable at each point $p\in \frz T^\ast\cup Z$.
In fact, if $\frz T^\ast\cup Z$ is normal crossing in $p$, since there are not
type changing components, $\tau_k$ is locally completable at $p$. But, in
general $\frz T^\ast\cup Z$ is not normal crossing. If $p_0$ is not normal
crossing we have, by construction 2.7, two irreducible components $H,H_1$ of
$\frz T^\ast$ and one irreducible component $H'$ of $Z$, with $[H']=[H]$,
meeting pairwise transversally at $p_0$.
\begin{center}\setlength{\unitlength}{1mm}\begin{picture}(40,40)
\put(0,20){\line(1,0){37}}
\put(20,2){\line(0,1){35}}
\put(35,5){\line(-1,1){30}}
\put(40,20){\makebox(0,0){$ _{H_1}$}}
\put(20,40){\makebox(0,0){$ _H$}}
\put(38,2){\makebox(0,0){$ _{H'}$}}
\put(10,10){\makebox(0,0){$+$}}
\put(30,30){\makebox(0,0){$+$}}
\put(30,30){\makebox(0,0){$+$}}
\put(30,15){\makebox(0,0){$-$}}
\put(10,25){\makebox(0,0){$-$}}
\put(20,-3){\makebox(0,0){\small Figure 2.8}}
\end{picture}
\end{center}
\vskip .2cm

\noindent Then by construction we have two signs $+1$ between $H$ and $H_1$
near $p_0$
(if not $H'$ would not cross $H$) and at most two signs $-1$ between $H$ and
$H'$ or between $H'$ and $H_1$ (see figure 2.8). So $\tau_k$ is locally
completable
at $p_0$.\\
\indent (3) If $H$ is a change component for $\tau_k$, then $H\subset
\frz T^\ast$ and if $[H]\not= 0$ then for some irreducible
component $Z_H$ of $Z$, $[H\cup Z_H]=0$ and  $Z_H\cap T^\ast=\emptyset$,
$Z_H\cap B_k=\emptyset$, by construction.\\
\indent So $\tau_k$ is completable for $k=1,...,m$. Let $P_k$ be a regular
function
inducing $\tau_k$. Then
$$T^\ast\subset\{P_1>0,\dots,P_m>0\},$$
\noindent but if $x\not\in T^\ast$, $x\in (\bigcup_{k=1}^mB_k)\cup\frz
T^\ast\cup Z$,
hence at least one among $P_k$ verifies $P_k(x)\leq 0$. So
$$T^\ast=\{P_1>0,\dots,P_m>0\}$$
\noindent then it is basic, hence $S$ is generically basic, but $S\cap\frz
S=\emptyset$
so, by 1.4, $S$ is basic.
\end{pf}

 From 2.8 we find for surfaces a geometric proof of a general basicness
characterization
in [AR2]. We call {\em birational model} of a semialgebraic set $S$ any
semialgebraic
set obtained from $S$ by a birational morphism on $X$.

\begin{cor}
Let $X$ be a surface and $S\subset X$ a semialgebraic set. Then, $S$ is
basic open if and only if $\frz S\cap S=\emptyset$ and for each birational
model $T$ of $S$ we have $\frz T^\ast\cap T^\ast$ is a finite set.
\end{cor}

\begin{pf}
By 1.4 we are done the {\em only if part}. Suppose now that $\frz S\cap
S=\emptyset$
and for each birational model $T$ of $S$ we have that $\frz T^\ast\cap T^\ast$
is a finite set.
Then, $S$ is open and we will prove that it is basic open.\\
\indent Take a compactification of $X$ (for instance its closure in a
projective
space) and then take a non-singular birational model $X_1$ of $X$, obtained by
a finite sequence of blowings-up along smooth centers. The strict
transform $S_1$ of $S$ is a birational model of $S$.\\
\indent Consider now the standard resolution $\pi:X_2\to X_1$ of $\frz S_1$
and take $S_2$ the strict transform of $S_1$ by $\pi$. Then $\frz S_2$ is
normal
crossing and $\frz S_2^\ast\cap S_2^\ast$ is a finite set, because $S_2$ is a
birational model of $S$. But, $\frz S_2^\ast\subset \frz S_2$ has not isolated
points (it is normal crossing) and $S_2^\ast$ is open, then $\frz S_2^\ast\cap
S_2^\ast=\emptyset$.
So by 2.5 and 2.8, $S_2^\ast$ is basic open. Hence $S$ is generically basic
and, as $\frz S\cap S=\emptyset$, $S$ is basic open.
\end{pf}

\section{Geometric review of fans}

For all notions of real algebra, real spectra, specialization, real
valuation rings, etc., we refer to [BCR]. Only for tilde operation we use a
slightly different definition: for a semialgebraic set $S$ in an algebraic
set $X$, $\tilde S$ is the constructible set of $\spr (\r(X))$ (instead of
${\cal P}(X)$) defined by the same formula which defines $S$. The properties
of this tilde operation are the the same as the usual ones (see [BCR, chap.7]).

Let $K$ be a real field, a subset $F=\{\al_1,\al_2,\al_3,\al_4\}$ of $\spr K$
is a {\em 4-element fan} (or simply a {\em fan}) if each $\al_i$ is the product
of the other three,
that is for each $f\in K$ we have
$$\al_i\al_j\al_k(f)=\al_l(f)$$
\noindent for all $\{i,j,k,l\}=\{1,2,3,4\}$, where $\al(f)$ denotes the sign
($1$ or $-1$)
of $f$ in the ordering $\al\in \spr(K)$.\\
\indent Given a fan $F$ we can find a valuation ring $V$ of $K$ such that\\
\indent a) Each $\al_i\in F$ is compatible with $V$; that is, the maximal
ideal $\gtm_V$ of $V$ is $\al_i$-convex.\\
\indent b) $F$ induces at most two orderings in the residue field $k_V$ of
$V$.\\
In this situation we say that $F$ trivializes along $V$ (see [BCR, chap.10] and
[Br]).

Let $X$ be a real algebraic set, and $\K(X)$ be the function field of $X$,
that is a finitely generated real extension of $\R$. Denote $K=\K(X)$.

\begin{defn}
A fan $F$ of $K$ is {\em associated to a real prime divisor $V$} if\\
\indent {\em (a)} $V$ is a discrete valuation ring such that $F$ trivializes
along $V$.\\
\indent {\em (b)} The residue field $k_V$ of $V$ is a finitely generated real
extension of $\R$ such that ${\rm dg.tr.}[K:\R]={\rm dg.tr.}[k_V:\R]+1$.
\end{defn}

\begin{rmk}
{\rm Let $F$ be a not trivial fan (i.e. the $\al_i$'s are distincts) associated
to a real prime divisor $V$, then it induces two distinct orderings
$\tau_1,\tau_2$
in $k_V$ ([BCR, 10.1.10]). If $F=\{\al_1,\al_2,\al_3,\al_4\}$ we suppose that
$\al_1,\al_3$ (resp. $\al_2,\al_4$) induce $\tau_1$ (resp. $\tau_2$) and we
write this
\[\begin{array}{ccccccccc}
V& &\al_1& &\al_3& &\al_2& &\al_4\\
\downarrow& & &\searrow\swarrow& & & &\searrow\swarrow& \\
k_V& & &\tau_1& & & &\tau_2&
\end{array}\]
Conversely, let $\tau_1,\tau_2\in\spr(k_V)$ be distinct, and let $t\in V$ be a
uniformizer for $V$. Each $f\in V$ can be written as $f=t^nu$, where $n$
is the valuation of $f$ and $u$ is a unit in $V$. Denote by $\ol u$ the class
of
$u$ in $k_V$ and consider the orderings in $K$ defined as follows:
\[\begin{array}{ccc}
\al_1(f)=\tau_1(\ol u)&;&\al_3(f)=(-1)^n\tau_1(\ol u)\\
\al_2(f)=\tau_2(\ol u)&;&\al_4(f)=(-1)^n\tau_2(\ol u)\\
\end{array}\]
\noindent They form a fan $F$ of $K$ associated to the real prime divisor
$V$.\\
\indent We may consider $\tau_1,\tau_2\in\spr(k_V)$ as elements of $\spr(V)$
with $\gtm_V$
as support. Then we have that $\al_1,\al_3$ (resp. $\al_2,\al_4$) specialize to
$\tau_1$ (resp. $\tau_2$) in $\spr(V)$.\\
\indent When $\al$ specializes to $\tau$, we write $\alpha\to \tau$.}
\end{rmk}

 From now we consider the field of rational functions $\K(X)$ of a compact
non-singular real
algebraic surface $X$, which is a finitely generated real extension of $\R$
with
transcendence degree over $\R$ equal to 2.

\begin{rmk}
{\rm Let $F$ be a fan in $\K(X)$ associated to a real prime divisor V of
$\K(X)$.
Then, $\r(X)\subset V$ (because $X$ is compact); consider the real prime
ideal $\gtp=\r(X)\cap\gtm_V$. We have that $V$ dominates $\r(X)_\gtp$ and there
are two possibilities:\\
\indent 1) If the height of $\gtp$ is 1, it is the ideal of an irreducible
algebraic curve $H\subset X$. Since $X$ is non-singular, $\r(X)_\gtp$ is
a discrete valuation ring, which is dominated by $V$. Hence $V=\r(X)_\gtp$
and $k_V$ is the function field $\K(H)$ of $H$; so $\tau_1,\tau_2\in\spr(H)$.\\
\indent 2) If $\gtp$ is a maximal ideal, it is the ideal of a point $p\in X$,
because $X$ is compact.}
\end{rmk}

\begin{defn}
Let $F$ be a fan of $\K(X)$ associated to a real prime divisor $V$ and let
$\gtp=\r(X)\cap\gtm_V$. The {\em center of $F$} is the zero set $Z(\gtp)$ of
$\gtp$.\\
\indent We say that {\em $F$ is centered at a curve (resp. a point)} if $\gtp$
has height 1
(resp. is maximal).
\end{defn}

\begin{lem}
Let $S\subset X$ be an open semialgebraic set. Then the following facts are
equivalent:\\ \indent {\em (i)} For each fan $F$ of
$\K(X)$ centered at a curve,  $\#(F\cap\tilde S)\not= 3$\\
\indent {\em (ii)} $\frz S$ has not positive type changing components with
respect
to the sign distributions $\s_i^S$ for $i=1,...,t$, defined in 2.4.
\end{lem}

\begin{pf}
Suppose to have a fan $F$ centered at a curve $H\subset X$, such that
$\#(F\cap\tilde S)=3$.
Then by remarks 3.2 and 3.3 we have:\\
\indent a) $F$ is associated to a real prime divisor $V=\r(X)_\gtp$, where
$\gtp$ is the ideal of $H$.\\
\indent b) If $F=\{\al_1,\al_2,\al_3,\al_4\}$, then
$\al_1,\al_3\to\tau_1\; ,\;\al_2,\al_4\to\tau_2$
in $\spr(V)$, with $\tau_1\not=\tau_2$ and $\tau_1,\tau_2\in\spr(\K(H))$.\\
\indent Suppose $\al_1,\al_2,\al_3\in \tilde S$ and $\al_4\not\in\tilde S$.
Remark
that an element of $\spr(\r(X)_\gtp)$ is a prime cone of $\spr(\r(X))$ which
support is contained in $\gtp$. So we can consider $\al_i,\tau_j\in\spr(\r(X))$
($i=1,2,3,4$, $j=1,2$) with $\tau_1,\tau_2\in\tilde H$ and
$\al_1,\al_3\to\tau_1\; ,\;\al_2,\al_4\to\tau_2$
in $\spr(\r(X))$.\\
\indent We have $\tau_1\in\ol{\tilde S}=\tilde{\ol S}$, because $\al_1,\al_3\in
\tilde S$.
But by [BCR, 10.2.8] there are precisely two prime cone different from $\tau_1$
in $\spr(\r(X))$
specializing to $\tau_1$, so they are $\al_1,\al_3$. And as
$\al_1,\al_3,\tau_1\in\tilde{\ol S}$,
we get that $\tau_1$ is an interior point of $\tilde{\ol S}$, so
$\tau_1\in\tilde{S^\ast}$.
This means $\tau_1\in\tilde H\cap\tilde{S^\ast}$, so $\di(H\cap S^\ast)= 1$.\\
\indent Now $\tau_2\in\tilde{\ol S}$, because $\al_2\in\tilde S$ and
$\al_2\to\tau_2$.
Again [BCR, 10.2.8] the prime cones specializing to $\tau_2$ and different
from it are precisely $\al_2,\al_4$. Since $\al_4\not\in\tilde S$ and
$\tilde S\cap\spr(\K(X))=\tilde{\ol S}\cap\spr(\K(X))$, we have that
$\al_4\not\in\tilde{\ol S}$, so $\tau_2$ is not interior to $\tilde{\ol S}$,
that
is $\tau_2\not\in\tilde{S^\ast}$. But $\ol S=\ol{S^\ast}$, then
$$\tau_2\in\widetilde{\ol{S^\ast}\setminus S^\ast}\subset\widetilde{\frz
S^\ast} \,.$$
\noindent This implies that $\di(H\cap\frz S^\ast)=1$, then $H$ is an
irreducible
component of $\frz S^\ast$. So $\di(S^\ast\cap\frz S^\ast)=1$ and by 2.5
$\frz S$ has a positive type changing component with respect $\s_i^S$ for some
$i=1,...,t$.\\
\indent Conversely, let $H$ be a irreducible component of $\frz S$ which is
positive type changing with respect $\s_i^S$ for some $i$. Then we can find
open sets $\om_1,\om_2\in H\cap{\rm Reg}(\frz S)$ such that\\
\indent a) $\om_1\subset\ol{\s_i^{-1}(1)}\cap\ol{\s_i^{-1}(-1)}=\ol S\cap\ol
A_i$\\
\indent b) $\om_2\subset{\rm Int}(\ol{\s_i^{-1}(1)})=S^\ast$\\
Let $\gtp$ be the ideal of $H$ in $\r(X)$ and $V$ be the discrete valuation
ring $\r(X)_\gtp$. Consider two orderings $\tau_1,\tau_2\in\spr(\K(H))$, with
$\tau_1\in \tilde\om_1$, $\tau_2\in\tilde\om_2$, and let $F$ be the fan defined
by $\tau_1,\tau_2$ as in 3.2. So
$\al_1,\al_3\to\tau_1\; ,\;\al_2,\al_4\to\tau_2$
in $\spr(\r(X))$, as before.\\
\indent We have $\al_2,\al_4\in\tilde S^\ast$, because $\tau_2\in\tilde S^\ast$
and $S^\ast$ is open. But $\al_2,\al_4\in\spr(\K(X))$ and
$\tilde S\cap\spr(\K(X))=\tilde S^\ast\cap\spr(\K(X))$, then
$\al_2,\al_4\in\tilde S$.\\
\indent On the other hand, $\tau_1\in \tilde{\ol S}\cap\tilde{\ol A_i}$. So
there
exist $\al\in\tilde S$ and $\beta\in\tilde A_i$ with $\al,\beta\to\tau_1$.
Again by [BCR, 10.2.8] we must have $\al=\al_1$ and $\beta=\al_3$. So
$\#(F\cap\tilde S)=3$.
\end{pf}

\begin{lem}
Let $S$ be a open semialgebraic set in $X$ such that $\frz S^\ast\cap S^\ast$
is
a finite set. Fix $p\in\partial S$. Then the following facts are equivalent:\\
\indent {\em (i)} For each fan $F$ centered at $p$,
$\#(F\cap\tilde S)\not= 3$.\\
\indent {\em (ii)} For each contraction
$\pi:X'\to X$ of a curve $E$ to the point $p$, no irreducible component of $E$
is positive type changing with
respect to $\s_i'=\s_i^S\cdot \pi$, for $i=1,...,t$.
\end{lem}

\begin{pf}
Suppose that there exist a contraction $\pi:X'\to X$ of a curve $E$ to the
point $p$
and $i=1,...,t$ such that an irreducible component $H$ of
$E$ is positive type changing with respect to $\s_i'$. But if $T=\pi^{-1}(S)$,
then $(\s_i')^{-1}(1)=T$, $(\s_i')^{-1}(-1)=\pi^{-1}(A_i)$; moreover, as
$p\not\in S$, the set $\{\pi^{-1}(A_i):i=1,...,t\}$ is precisely the set
of connected component of $X'\setminus(T\cup\frz T)$, then $\s_i'=\s_i^T$.
Now by 3.5 there exists a fan $F=\{\al_1,\al_2,\al_3,\al_4\}$ of $\K(X')$ with
center the curve $H$ such that $\#(F\cap\tilde T)=3$. \\
\indent The contraction $\pi$ induces a field isomorphism
$$\pi_*:\K(X)\to\K(X')$$
\noindent and an injective ring homomorphism
$\pi_*\,\vline_{\r(X)}:\r(X)\to\r(X').$
Let $G$ be the fan of $\K(X')$ inverse image of $F$ by $\pi_*$, namely
$$G=\{\pi_*^{-1}(\al_1),\pi_*^{-1}(\al_2),\pi_*^{-1}(\al_3),\pi_*^{-1}(\al_4)\}\subset\spr(\K(X))$$
Then $\#(G\cap\tilde S)=3$ and we have to prove that $p$ is the center
of $G$. Let $V$ be the real prime divisor associated to $F$, then
$\pi_*^{-1}(V)=W$
is the real prime divisor associated of $G$, so
$$\gtm_W\cap\r(X)=\pi_*^{-1}({\cal J}(H))$$
\noindent where ${\cal J}(H)$ denotes the ideal of $H$ in $\r(X')$. Hence,
$$\gtm_W\cap\r(X)={\cal J}(\pi(H))=\gtm_p$$
\noindent with $\gtm$ the maximal ideal of $p$.\\
\indent On the contrary, we suppose that no irreducible component of $E$ is
positive type changing with respect to $\s_i'=\s_i\cdot\pi$, for each
contraction $\pi$ of a curve to $p$. Take a neighbourhood $U$ of $p$,
homeomorphic to a disk in $\R^2$ ($X$ is non-singular), such that $U$ does not
meet any irreducible component of $\frz S$ unless it contains $p$. Consider the
sign distributions $\delta_j$ in $X\setminus (\frz S\cup \partial U)$ for
$j=1,...,l$, defined by
\begin{eqnarray*}
\delta_j^{-1}(1)&=&U\cap S\\
\delta_j^{-1}(-1)&=&B_j
\end{eqnarray*}
\noindent where $B_1,\dots,B_l$ are the connected components of
$U\setminus(S\cup\frz S)$.
As $\frz S^\ast\cap S^\ast$ is a finite set, by 2.5 $\frz S$ has not positive
type changing components with respect to $\s_i^S$ for all $i=1,...,t$.
We claim that $\frz(U\cap S)$ has no type changing components at all. In fact,
as $B_j\subset A_i$ for some $i$, there are not positive ones; $\partial U$
cannot be type changing because the signs may lie only on one side of it and
no other component can be negative type changing, for the same reasons as in
the proof of 2.6.\\
\indent Now if $\frz S$ is normal crossing at $p$, by 2.2 $\delta_j$ is locally
completable at $p$ for $i=1,...,l$. If not, consider the standard singularity
resolution $\pi:X'\to X$ of $\frz S$ at $p$ and the sign distributions
\[\begin{array}{c}
\s_i'=\s_i^S\cdot\pi,\;{\it for}\;i=1,...,t\\
\delta_j'=\delta_j\cdot\pi,\;{\it for}\;i=1,...,l
\end{array}   \]
\noindent As no irreducible component of $\pi^{-1}(p)$ is positive type
changing
with respect to $\s_i'$ for all $i$, the same is true with respect to
$\delta_j'$ ($j=1,...,l$). More over each such component has a M\"obius
neighbourhood
in $\pi^{-1}(U)$ where the sign minus can occur locally only on one side of the
curve, so it cannot be negative type changing.
Then, by 2.2 $\delta_j$ is locally completable at $p$ for each $j=1,..,l$.\\
\indent Hence, for each $j=1,...,l$, take $f_j\in\r(X)$ and an open set $U_j\ni
p$, $U_j\subset U$, such that $f_j$ induces
$\delta_j$ on $U_j$. Consider $A=\bigcap_{j=1}^lU_j$, then
$$S\cap A=\{f_1>0,\dots,f_l>0\}\cap A.$$
\noindent In fact, by definition of locally completable we have
$$S\cap A\subset \{f_1>0,\dots,f_l>0\};$$
\noindent but if $x\in A\setminus S$, then $x\in(\bigcup B_j)\cup\frz S$, so
there
is a $j_0$ such that $f_{j_0}(x)\leq 0$.\\
\indent Let $F=\{\al_1,\al_2,\al_3,\al_4\}$ a fan centered at $p$. Clearly
$F\subset\tilde A$, because each $\al_i$ specializes in $\spr(\r(X))$ to the
prime  cone having $\gtm_p$ as support, i.e the unique prime cone giving to
$f\in\r(X)$ the sign of $f(p)$ ([BCR, 10.2.3]). Suppose that
$\al_1,\al_2,\al_3\in\tilde S$
and $\al_4\not\in\tilde S$; then $\al_1,\al_2,\al_3\in\tilde A\cap\tilde S$
and $\al_4\not\in\tilde A\cap\tilde S$. So for all $j=1,...,l$, $\al_i(f_j)>0$
for $i=1,2,3$ and there is $j_0$ such that $\al_4(f_{j_0})<0$. But this is
imposible because $F$ is a fan and
$\al_1\al_2\al_3(f_{j_0})\not=\al_4(f_{j_0})$.
\end {pf}

\begin{rmk}
{\rm Let $S$ be an open semialgebraic set such that $\frz S^\ast\cap S^\ast$ is
a
finite set. Let $p\in\frz S$ such that $p\not\in\partial S$ or $\frz S$ is
normal crossing at $p$. Then for each fan centered at $p$,
$\#(F\cap\tilde S)\not= 3$.}
\end{rmk}

\begin{thm} {\rm (See [Br] and [AR1])}
Let $X$ be a real irreducible algebraic surface. Let $S$ be a semialgebraic set
such that $\frz S\cap S=\emptyset$ (resp. a closed semialgebraic set).
Then, $S$ is basic open (resp. basic closed) if and only if for each fan $F$ of
$\K(X)$ which is associated to a real prime divisor, $\#(F\cap\tilde S)\not=3$.
\end{thm}

\begin{pf}
Suppose $S$ to be basic open, then
$$S=\{f_1>0,\dots,f_r>0\},\;{\rm with}\; f_1,\dots,f_r\in\r(X).$$
\noindent Let $F=\{\al_1,\al_2,\al_3,\al_4\}$ be a fan in $\K(X)$ and suppose
that $\al_1,\al_2,\al_3\in S$ and $\al_4\not\in S$. Then for all $i=1,...,r$,
$f_i(\al_j)>0$ ($j=1,2,3$) and there exists $i_0\in\{1,...,r\}$ such that
$\al_4(f_{i_0})>0$; which is impossible because $F$ is a fan.\\
\indent Conversely, suppose $S$ to be not basic open. If $X$ is compact and non
singular, by 2.8, 3.5 and 3.6 we have done. If not, take a birational model
$X_1$ of $X$ obtained compactifying and desingularizing $X$. Let $S_1$ be
the strict transform of $S$ in $X_1$. Then $S_1$ is not basic open and
$\frz(S_1)\cap S_1=\emptyset$,
because $S$ verifies these properties. So, by
2.8, 3.5 and 3.6, we can find a fan $F$ of $\K(X_1)$ associated to a real prime
divisor
such that $\#(F\cap\tilde S_1)=3$.
Since $\K(X)$ and $\K(X_1)$ are isomorphic, $F$ gives a fan $G$ of $\K(X)$ such
that $G$ is associated to a real prime divisor and $\#(G\cap\tilde S)=3$.
\end{pf}

\section{The algorithms}

By 2.3 (see also [ABF]) there is an algorithmic method for checking properties
{\em (a)} and
{\em (b)} of 2.8; so we can decide algorithmically if a semialgebraic $S$ with
$\frz S\cap S=\emptyset$ is open basic. This method works as follows:\\
\indent 1) It calculates $\di(S^\ast\cap\frz S^\ast)$ by tecniques of
cilindrical algebraic descomposition (C.A.D.), for instance (see [BCR]). If it
is equal to 1, we
know that $S$ is not basic. If not, we continue with 2).\\
\indent 2) It decides if some irreducible component of the exceptional divisor
of the standard resolution of $\frz S$ is positive type changing using
\framebox{A1} and \framebox{A2}
in the points of $\partial S$ which are not normal crossing (remark that
\framebox{A2}
decides if the non-local completability at a point is due to a positive or a
negative type changing component after some blowing-up).

Moreover, from [Vz] we have a complete description for fans in $\R(x,y)$
associated to
a real prime divisor.\\
\indent Consider the field  $\R((u,v))$ of formal series in two variables over
$\R$, with the ordering which
extends $0^+$ in $\R((u))$ by $v>0$. Any ordering in $\R(x,y)$  is defined
by an ordered $\R$-homomorphism $\psi:\R(x,y)\to\R((u,v))$ (see [AGR]). So a
non-trivial fan $F$
is given by 4 homomorphisms $\psi_1,\psi_2,\psi_3,\psi_4$. More precisely:

\begin{thm}
{\rm (See [Vz])} Let $F$ be a fan in $\R(x,y)$. Then $F$ is described as
follows:\\
\indent {\em 1)} If $F$ has as center an irreducible curve $H\subset{\Bbb S}^2$
and if $P(x,y)\in\R[x,y]$
is a polynomial generating the ideal ${\cal J}(H)\subset\R[x,y]$ of the image
of $H$ by a suitable stereographic projection, then
$$\psi_i:\R(x,y)\to \R((t,z)),\; for\; i=1,2,3,4$$
\noindent are defined (possibly interchanging $x$ and $y$) as follows:
\[\left\{ \begin{array}{l}
\psi_1(x)=a_1+\delta t^N\\ \psi_1(y)=a_2+\sum_{i\geq 1}c_it^{n_i}+z
\end{array}\right. \qquad \left\{\begin{array}{l}
\psi_2(x)=b_1+\delta't^M\\ \psi_2(y)=b_2+\sum_{i\geq 1}d_it^{m_i}+z
\end{array}\right.\]
\[\left\{ \begin{array}{l}
\psi_3(x)=a_1+\delta t^N\\ \psi_3(y)=a_2+\sum_{i\geq 1}c_it^{n_i}-z
\end{array}\right. \qquad \left\{\begin{array}{l}
\psi_4(x)=b_1+\delta't^M\\ \psi_4(y)=b_2+\sum_{i\geq 1}d_it^{m_i}-z
\end{array}\right.\]
\noindent where $(a_1+\delta t^N,\, a_2+\sum_{i\geq 1}c_it^{n_i})$ and
$(b_1+\delta't^M,\, b_2+\sum_{i\geq 1}d_it^{m_i})$ are irreducible Puiseux
parametrizations
of two half-branches of $H$, centered respectively at $(a_1,a_2)$, $(b_1,b_2)$.

{\em 2)} If $F$ is centered at a point $p\in{\Bbb S}^2$, we may suppose
$p=(0,0)$
in a suitable stereographic projection, then
$$\psi_i:\R(x,y)\to \R((z,t)),\; for\; i=1,2,3,4$$
\noindent are given by one of the following expressions:\\
\indent {\em a)}
\[\left\{ \begin{array}{l}
\psi_1(x)=t\\ \psi_1(y)=tw \phantom{-}
\end{array}\right. \qquad \left\{\begin{array}{l}
\psi_2(x)=tw'\phantom{-}\\ \psi_2(y)=t
\end{array}\right.\]
\[\left\{ \begin{array}{l}
\psi_3(x)=-t\\ \psi_3(y)=-tw
\end{array}\right. \qquad \left\{\begin{array}{l}
\psi_4(x)=-tw'\\ \psi_4(y)=-t
\end{array}\right.\]
\noindent with $w\in\{z+a,-z+a:a\in\R\}$, $w'\in\{z,-z\}$.\\
\indent {\em b)} Up to interchanging $x$ and $y$,
\[\left\{ \begin{array}{l}
\psi_1(x)=\delta t^N\\
\psi_1(y)=\sum_{i=1}^sc_it^{n_i}+t^mw\phantom{(-1)^{n_i}(-1)^m}
\end{array}\right. \qquad \left\{\begin{array}{l}
\psi_2(x)=\delta t^N\\
\psi_2(y)=\sum_{i=1}^sc_it^{n_i}+t^mw'\phantom{(-1)^{n_i}(-1)^m}
\end{array}\right.\]
\[\left\{ \begin{array}{l}
\psi_3(x)=(-1)^N\delta t^N\\
\psi_3(y)=\sum_{i=1}^s(-1)^{n_i}c_it^{n_i}+(-1)^mt^mw
\end{array}\right. \qquad \left\{\begin{array}{l}
\psi_4(x)=(-1)^N\delta t^N\\
\psi_4(y)=\sum_{i=1}^s(-1)^{n_i}c_it^{n_i}+(-1)^mt^mw'
\end{array}\right.\]
\noindent with $\delta\in\{1,-1\}$; $c_i\in\R$ for $i=1,...,s$; $N\leq
n_1<n_2<\dots <n_s$,
${\rm g.c.d.}(N,n_1,...,n_s)=d$ and ${\rm g.c.d}(d,m)=1$; and
$w,w'\in\{z+a,-z+a,1/z,-1/z:a\in \R\}$, with $w\not= w'$ if $d$ is odd and
$w\not= w',w\not= -w'$ if $d$ is even. If $N=1$, $c_1=\dots=c_s=0$, then
$w,w'\not\in\{1/z,-1/z\}$.\\
\indent {\em c)} Up to interchanging $x$ and $y$,
\[\left\{ \begin{array}{l}
\psi_1(x)=\delta t^N\\ \psi_1(y)=\sum_{i=1}^sc_it^{n_i}+t^mw
\end{array}\right. \qquad \left\{\begin{array}{l}
\psi_2(x)=(-1)^{N/d}\delta t^N\\
\psi_2(y)=\sum_{i=1}^s(-1)^{n_i/d}c_it^{n_i}+t^mw'
\end{array}\right.\]
\[\left\{ \begin{array}{l}
\psi_3(x)=\delta t^N\\ \psi_3(y)=\sum_{i=1}^sc_it^{n_i}-t^mw
\end{array}\right. \qquad \left\{\begin{array}{l}
\psi_4(x)=(-1)^{N/d}\delta t^N\\
\psi_4(y)=\sum_{i=1}^s(-1)^{n_i/d}c_it^{n_i}-t^mw'
\end{array}\right.\]
\noindent with $\delta,c_i,N,n_i,w,w'$ as in {\em b)}, but $d$ always even and
without supplementary conditions on $w,w'$.
\end{thm}

\begin{rmk}
{\rm To each fan $F$ in $\R(x,y)$ centered at $(0,0)$ we can associate two
families of arcs through $(0,0)$ which are parametrized by $z\in (0,\epsilon)$.
In fact, for each fixed $z\in(0,\epsilon)$, $\psi_1$ and $\psi_3$ (resp.
$\psi_2$
and $\psi_4$) define the two half-branches of the same curve germ $\gamma_1^z$
(resp. $\gamma_2^z$). This curves, $\gamma_1^z$ and $\gamma_2^z$, verify one
with respect to the other conditions a) and b) of 2.3 (see [ABF, 2.19]).}
\end{rmk}

We want to know the relations between these two families of arcs associated to
$F$ and the families of arcs of \framebox{A1} and \framebox{A2} (2.3).

Let $C$ be a curve germ through $(0,0)\in\R^2$; and consider the standard
resolution
$\pi=\pi_N\cdots\pi_1$ of $C$ at $(0,0)$. For each $i=1,\dots,N$, denote by
$C_i$ the curve $(\pi_i\cdots\pi_1)^{-1}(C)$, by $D_i$
the exceptional divisor arising during the $i^{th}$ blowing-up, and by
$E_i$ the exceptional curve after $i$ blowings-up (i.e.
$E_1=(\pi_i\cdots\pi_1)^{-1}(0,0)$). $D_i$ is an irreducible component of
$E_i$.

\begin{defn}
Let $F$ be a fan of $\R(x,y)$ with center $(0,0)\in\R^2$, for $i=1,...,N$
denote by $F_i$ the fan obtained from $F$ after $i$ blowings-up by lifting the
orderings of $F$. We say that $F$ has the {\em property $\star(\rho)$} with
respect to $C$ if it verifies:\\
\indent {\em a)} $F_{\rho-1}$ is centered at the point $0=C_{\rho-1}\cap
D_{\rho-1}$.\\
\indent{\em b)} $F_{\rho-1}$ is described in the sense of 4.1 as in {\em 2-a)}
or {\em 2-b)}
with $N=1$, $c_1=\dots=c_s=0$.\\
\indent {\em c)} $C_{\rho-1}$ is not tangent to any curve of the two families
associates to $F_{\rho-1}$.
\end{defn}

\begin{rmk}
{\rm A fan $F$ verifies $\star(\rho)$ with respect to a curve $C$ if and only
if for each
$z\in(0,\epsilon)$ the curves $\gamma_1^z$ and $\gamma_2^z$ verify a) and b) of
2.3 with respect to $C$ (see
also [ABF, 2.9]).}
\end{rmk}

Finally we have:

\begin{thm}
Let $S$ be an open semialgebraic set in $X$, such that $\frz S\cap
S=\emptyset$. Suppose that $\frz S^\ast\cap S^\ast$
is a finite set and $S$ is not basic. Then, there exists an algorithmic
method for finding a fan $F$ of $\K(X)$ with
$\#(F\cap\tilde S)=3$.
\end{thm}

\begin{pf}
By 2.5 $\frz S$ has not type changing components with respect to $\s_i^S$; then
by 2.8 there is at least one irreducible component $D_\rho$ of the exceptional
divisor
of a standard resolution $\pi$ of $\frz S$ at a point $O=(0,0)$ which is
positive type changing
with respect to some $\s_i'=\s_i^S\cdot\pi$. By 2.3 we can find algorithmically
two arcs $\gamma_1,\gamma_2$ with the properties a) and b) of 2.3 with respect
to an irreducible component of $\frz S$ through $p=\pi(D)$, for $\rho>0$, such
that $\gamma_1$
joins two regions in $(\s_i^S)^{-1}(1)$, while $\gamma_2$ joins a region in
$(\s_i^S)^{-1}(1)$ to a region in $(\s_i^S)^{-1}(-1)$. Moreover, each
$\gamma_i$
($i=1,2$) is defined by open conditions.\\
\indent Then $(\gamma_1)_\rho$ and $(\gamma_2)_\rho$ are smooth arcs wich meet
$D_\rho$ transversally in different points $p,q\in D_\rho$ and $\rho$ is the
first level in the resolution process at wich $\gamma_1$ and $\gamma_2$ are
separated.
By [ABF, 2.19], $\gamma_1$ and $\gamma_2$ are parametrized by
\[ \gamma_1:\left\{ \begin{array}{l}
x=\delta t^N\\ y=\sum_{i=1}^sc_it^{n_i}+f(t)
\end{array}\right. \qquad \gamma_2:\left\{\begin{array}{l}
x=\delta' t^N\\ y=\sum_{i=1}^sd_it^{n_i}+g(t)
\end{array}\right.\]
\noindent where $\delta,\delta'\in\{1,-1\}$, $c_i,d_i\in\R$ are determined by
 [ABF, 2.19] as follows: $\delta'=\delta$ and $d_i=c_i$  or
$\delta'=(-1)^{N/d}\delta$ and $d_i=(-1)^{n_i/d}c_i$, for $i=1,...,s$, with
$d={\rm g.c.d.}(N,n_1,\dots,n_s)$;
and $f(t)=at^m+...$, $g(t)=bt^m+...$ with $a\not= b$.\\
\indent We can construct four fans with the property $\star(\rho)$ with respect
to
$\gamma_1$ and $\gamma_2$ as follows:\\
\indent At the level $\rho$ of the resolution process we have four fans
centered
at $D_\rho$ in half-branches at $p$ and $q$ (Fig. 4.5), obtained by taking
respectively a half-branch
of $D_\rho$ at $p$ and an other at $q$.
\begin{center}\setlength{\unitlength}{1mm}\begin{picture}(100,30)
\multiput(20,2)(20,0){4}{\line(0,1){26}}
\multiput(20,30)(20,0){4}{\makebox(0,0){$_{D_\rho}$}}
\multiput(16,7)(0,1){3}{\line(4,1){4}}
\multiput(56,7)(0,1){3}{\line(4,1){4}}
\multiput(20,8)(0,1){3}{\line(4,-1){4}}
\multiput(60,8)(0,1){3}{\line(4,-1){4}}
\multiput(16,21)(0,1){3}{\line(4,-1){4}}
\multiput(36,21)(0,1){3}{\line(4,-1){4}}
\multiput(20,20)(0,1){3}{\line(4,1){4}}
\multiput(40,20)(0,1){3}{\line(4,1){4}}
\multiput(36,11)(0,1){3}{\line(4,-1){4}}
\multiput(76,11)(0,1){3}{\line(4,-1){4}}
\multiput(40,10)(0,1){3}{\line(4,1){4}}
\multiput(80,10)(0,1){3}{\line(4,1){4}}
\multiput(56,17)(0,1){3}{\line(4,1){4}}
\multiput(76,17)(0,1){3}{\line(4,1){4}}
\multiput(60,18)(0,1){3}{\line(4,-1){4}}
\multiput(80,18)(0,1){3}{\line(4,-1){4}}
\multiput(20,10)(20,0){4}{\circle*{1.5}}
\multiput(20,20)(20,0){4}{\circle*{1.5}}
\multiput(22,12)(40,0){2}{\makebox(0,0){$_p$}}
\multiput(22,18)(20,0){2}{\makebox(0,0){$_q$}}
\multiput(42,8)(40,0){2}{\makebox(0,0){$_p$}}
\multiput(62,22)(20,0){2}{\makebox(0,0){$_q$}}
\put(50,-1){\makebox(0,0){Figure 4.5}}
\end{picture}
\end{center}
\indent Going back by the birational morphism $\pi_\rho\cdots\pi_1$ we obtain
four fans
centered at $O=(0,0)$ such that all they have the property $\star(\rho)$ with
respect to $\gamma_1$ and $\gamma_2$. More precisely applying again
[ABF, 2.19] we can describe them in terms of 4.1: for every pair
$\eta,\eta'\in\{1,-1\}$
we have one of this fans $F_{\eta,\eta'}$ and his associated arcs are
\[ \gamma_1^z:\left\{ \begin{array}{l}
x=\delta t^N\\ y=\sum_{i=1}^sc_it^{n_i}+(\eta z+a)t^m
\end{array}\right. \qquad \gamma_2^z:\left\{\begin{array}{l}
x=\delta' t^N\\ y=\sum_{i=1}^sd_it^{n_i}+(\eta'z+b)t^m
\end{array}\right.\]
\indent By the construction of $\gamma_1,\gamma_2$ it is easy to check that
$\#(F_{\eta,\eta'}\cap\tilde S)=3$.
\end{pf}

\begin{rmk}
{\rm In the hypotesis of 3.13 there are in fact infinite fans $F$ of $\K(X)$
verifying
$\#(F\cap\tilde S)=3$, because there are infinite pairs of arcs joining
respectively two region with sign $1$ and a region with sign $1$ with a region
with sign $-1$.\\
\indent So we find only fans with $w,w'\in\{z+a,-z+a:a\in\R\}$, according to
description 4.1.2. In other case, $(\gamma_1^z)_{\rho}$ or
$(\gamma_2^z)_{\rho}$
for some $\rho$, would be tangent to $D_{\rho}$, but applying \framebox{A1} and
\framebox{A2} as in [ABF] we take $\gamma_1$, $\gamma_2$ without this property.
This means that if it exists a fan $F$ with $w\;{\rm or}\; w'\in\{1/z,-1/z\}$
and $\#(F\cap\tilde S)=3$, there is another fan $F'$ with
$w,w'\not\in\{1/z,-1/z\}$
and $\#(F'\cap\tilde S)=3$.}
\end{rmk}

\section{Principal sets}

 Using the results of the previous sections, we obtain a simple
characterization
of principal open (resp. closed) sets. In order to conserve the unity of this
paper we give all results about principal sets in dimension 2, but in fact
they can be extended to arbitrary dimension following similar proofs (Remarks
5.9).
Details can be found in [Vz].\\
\indent Let $X$ be a compact, non singular, real algebraic surface.

\begin{defn}
A semialgebraic set $S\subset X$ is {\em principal open} (resp. {\em principal
closed}) if there exists $f\in\r(X)$ such that
$$S=\{x\in X:f(x)>0\}$$
$$(resp.\; S=\{x\in X:f(x)\geq 0\})$$
\end{defn}

\begin{defn}
A semialgebraic set $S\in X$ is {\em generically principal} if there exists a
Zariski closed set $C\in X$ with $\di (C)\leq 1$ such that $S\setminus C$ is
principal open.
\end{defn}

\begin{rmk}
{\rm A semialgebraic set $S$ is principal closed if and only if $X\setminus S$
is principal open.\\
\indent Then it suffices to work with principal open sets.}
\end{rmk}

\begin{nott}
{\rm Let $S$ be an open semialgebraic set and $Y=\frz S$. We denote by $S^c$
the open
semialgebraic set $X\setminus(S\cup Y)$ and by $A_1,\dots,A_t$ (resp.
$B_1,\dots,B_l$) the connected component of $S^c$ (resp. $S\setminus Y$).\\
\indent Let $\s_i$ be the sign distributions $\s_i^S$ for $i=1,\dots,t$ and
$\s_j^c$ be the sign distributions $\s_j^{S^c}$ for $j=1,\dots,l$ defined as
in 2.4. And denote by $\delta$ the total sign distribution defined by
\begin{eqnarray*}
\delta^{-1}(1)&=&S\setminus Y\\
\delta^{-1}(-1)&=&S^c
\end{eqnarray*}               }
\end{nott}

\begin{rmk}
{\rm A semialgebraic set $S$ such that $\frz S\cap S=\emptyset$ is principal
open if and only if the sign distribution $\delta$ is admissible (that is,
there exists $f\in\r(X)$ such that $f$ induces $\delta$ on $X\setminus\frz
S$).}
\end{rmk}

\begin{thm}
Let $S$ be a semialgebraic set such that $\frz S\cap S=\emptyset$.
Then $S$ is principal open if and only if $S^*\cap\frz S^*$ and
$(S^c)^*\cap\frz(S^c)^*$ are finite sets.
\end{thm}

\begin{pf}
Suppose that $S$ is principal then $S$ and $S^c$ are basic and by 2.8 no
irreducible component of $\frz S$ is positive type changing with respect
to $\s_i$ and $\s_j^c$ for each $i=1,...,t$, $j=1,...,l$. Applying now 2.5
we have that $S^*\cap\frz S^*$ and $(S^c)^*\cap\frz(S^c)^*$ are finite sets.\\
\indent Conversely suppose that $S^*\cap\frz S^*$ and $(S^c)^*\cap\frz(S^c)^*$
are finite sets,
then by 2.5 again no irreducible component of $\frz S$ is positive type
changing
with respect to $\s_i$ and $\s_j^c$ for $i=1,...,t$, $j=1,...,l$.\\
\indent Remark that an irreducible component $H$ of $\frz S$ is positive (resp.
negative) type
changing with respect to $\delta$ if and only if $H$ is positive type changing
with respect to $\s_i$ for some $i$ (resp. $\s_j^c$ for some $j$).\\
\indent Hence no irreducible component of $\frz S$ is type changing with
respect
to $\delta$. Denote by $Z^c$ the union of all the change components of $\frz S$
and by $Z$ the set of points where $Z^c$ has dimension 1. So $[Z]=[Z^c]$ and
$[Z]=0$ in ${\rm H}_1(M,\Z_2)$, because it bounds the open sets ${\rm
Int}(\ol{\s^{-1}(1)})$
and ${\rm Int}(\ol{\s^{-1}(-1)})$. Then by [BCR, 12.4.6] the ideal ${\cal
J}(Z^c)$
of $Z^c$ is principal. Let $f$ be a generator of ${\cal J}(Z^c)$. Again by
[BCR, 12.4.6] for each irreducible component $H_k$ of $\frz S$ not lying in
$Z^c$ we can choose a generator $h_k$ of ${\cal J}(H_k)^2$ (wich exists because
$2[H_k]=0$). Then the regular function $f\cdot\prod h_k$ induces $\delta$
or $-\delta$. So $\delta$ is admisible and $S$ principal.
\end{pf}

Remark that this proof is almost the same as the proof of [AB, Proposition 2].

\begin{thm}
Let $S$ be a semialgebraic set in $X$.\\
\indent (1) $S$ is principal open if and only if $\frz S\cap S=\emptyset$ and
for each fan $F$ centered at a curve, $\#(F\cap\tilde S)\not= 1,3$.\\
\indent (2) $S$ is principal closed if and only if $\frz S\cap(X\setminus
S)=\emptyset$
and for each fan $F$ centered at a curve,  $\#(F\cap\tilde S)\not= 1,3$.
\end{thm}

\begin{pf}
It is immediately using 5.6 and 3.5.
\end{pf}

\begin{rmks}
{\rm (1) Results 5.6 and 5.7 can be generalized to an arbitrary surface
compactifying and desingularizing as in 2.9.\\

\indent (2) All this section can be generalized to a compact, non singular,
real
algebraic set $X$, because the results of the previus sections used here
(specifically
1.3 and 2.5) can be generalized to arbitrary dimension. Moreover, defining fan
 centered at a hypersurface $H$ of $X$ as a fan $F$ associated to a real prime
divisor $V$ such that the prime ideal $\gtp=\gtm_V\cap\r(X)$ has height $1$
and ${\cal Z}(\gtp)=H$, we find an improvement of 4-elements fans criterion
[Br, 5.3].\\

\indent (3) For a compact, non singular, real algebraic set we obtain:\\
{\em A semialgebraic set $S$ is principal open (resp. closed) if and only
if $\frz S\cap S=\emptyset$ (resp. $\frz S\cap(X\setminus S)=\emptyset$) and
$S$ is generically principal.}\\      }
\end{rmks}

\end{document}